\documentclass[12pt,a4paper]{amsart}
\usepackage[a4paper,left=2.2cm,right=2.2cm,top=2.0cm,
  bottom=2.0cm]{geometry}
\usepackage{amsfonts,amscd,amssymb,amsmath,amsthm,hyperref,mathrsfs,multirow,xcolor,lscape,longtable,pbox,lipsum}
\usepackage[thinlines]{easytable}
\usepackage{microtype}
\usepackage{mathtools}
\usepackage{eucal}
\usepackage[T1]{fontenc}
\usepackage{lmodern}
\usepackage[english]{babel}
\usepackage{microtype}
\newtheorem{theorem}{Theorem}[section]
\newtheorem{definition}[theorem]{Definition}
\newtheorem{lemma}[theorem]{Lemma}
\newtheorem{proposition}[theorem]{Proposition}
\newtheorem{corollary}[theorem]{Corollary}
\newtheorem{remark}[theorem]{Remark}

\newtheorem{problem}[theorem]{Problem}

 \usepackage[all,cmtip]{xy}\usepackage{tikz}\usetikzlibrary{arrows}
\usetikzlibrary{matrix}\usetikzlibrary{calc}

\newcommand{\ZZ}{\mathbb{Z}}

\newcommand{\FF}{\mathbb{F}}

\newcommand{\HH}{\mathbb{H}}
\newcommand{\QQ}{\mathbb{Q}}
\newcommand{\NN}{\mathbb{N}}
\newcommand{\RR}{\mathbb{R}}

\newcommand{\CC}{\mathbb{C}}

\newcommand{\rT}{\mathrm{T}}

\newcommand{\Per}{\mathrm{Per}}

\newcommand{\id}{\mathrm{id}}
\newcommand{\sff}{\mathrm{sf}}

\newcommand{\BL}{\mathrm{BL}}

\newcommand{\het}{\mathrm{ht}}

\newcommand{\Span}{\mathrm{Span}}
\newcommand{\disc}{\mathrm{disc}}
\newcommand{\Comm}{\mathrm{Comm}}

\newcommand{\cC}{\mathcal{C}}

\newcommand{\cG}{\mathcal{G}}

\newcommand{\cI}{\mathcal{I}}

\newcommand{\cK}{\mathcal{K}}

\newcommand{\sM}{\mathscr{M}}

\newcommand{\cO}{\mathcal{O}}
\newcommand{\cP}{\mathcal{P}}

\newcommand{\cR}{\mathcal{R}}
\newcommand{\cS}{\mathcal{S}}
\newcommand{\cT}{\mathcal{T}}

\newcommand{\ba}{\mathbf{a}}

\newcommand{\dd}{\mathrm{d}}

\newcommand{\GL}{\mathrm{GL}}
\newcommand{\Lip}{\mathrm{Lip}}
\newcommand{\tr}{\mathrm{tr}}

\newcommand{\Aut}{\mathrm{Aut}}
\newcommand{\SL}{\mathrm{SL}}
\newcommand{\Res}{\mathrm{Res}}
\newcommand{\diag}{\mathrm{diag}}
\newcommand{\supp}{\mathrm{supp}}
\newcommand{\Stab}{\mathrm{Stab}}
\newcommand{\per}{\mathrm{per}}

\newcommand{\red}[1]{\textcolor[rgb]{1.00,0.00,0.00}{#1}}

\newcommand{\blue}[1]{{\color{blue} {#1}}}

\title[Density of periodic measures with arithmetic restrictions]{Density of periodic measures on $\SL_2(\ZZ)\backslash\SL_2(\RR)$\\ with prescribed arithmetic information.}

\begin{document}

\author{Lifan Guan}
\address{Institute for Theoretical Sciences, School of Science, Westlake University,
600 Dunyu Road, Sandun town, Xihu district, Hangzhou 310030, Zhejiang Province, China.
}
\email{guanlifan@westlake.edu.cn}

\author{Chengyang Wu}
\address{Department of Mathematics, University of Chicago, Chicago, IL, US, 60637}
\email{chengyangwu1999@gmail.com}

\author{Weisheng Wu}
\address{School of Mathematical Sciences, Xiamen University, Xiamen, 361005, P.R. China}
\email{wuweisheng@xmu.edu.cn}

\date{9 November 2026}

%----------------------------------------------------------------------------------------
%	ABSTRACT
%----------------------------------------------------------------------------------------
\begin{abstract}
Let $F$ be the diagonal flow on the modular space $X_2$, and let $\cS_f$ be the set of positive square-free numbers, which is identified with the set of all real quadratic fields via the natural map $d\mapsto \QQ(\sqrt d)$. In this paper, we prove that, for any finite subset $A\subseteq\cS_f$, there exists a density-zero set $B\subseteq\cS_f\setminus A$ such that the set of all $F$-periodic measures with associated quadratic fields from $B$ is dense in the set of all $F$-invariant measures on $X_2$.
\end{abstract}

\maketitle

%\tableofcontents

%----------------------------------------------------------------------------------------
%	INTRODUCTION
%----------------------------------------------------------------------------------------

\section{Introduction}\label{S:intro}
The goal of this paper is to understand the density of periodic measures\footnote{Throughout this paper, all measures are assumed to be probability measures.} with arithmetic restrictions for the geodesic flow on modular surfaces. For hyperbolic flows on compact smooth manifolds, it is well known that (see \cite[Theorem 1]{Sigmund72}):

\begin{theorem}\label{T:Sigmund}
The set of all periodic measures is dense, under the weak-${}^*$ topology, in the set of all invariant measures for any transitive Anosov flow on any compact smooth manifold.
\end{theorem}
%\red{Shall we restrict to a compact smooth manifold? For non-uniform hyperbolic case, shall we restrict to a Pesin set? }

In 2010, Coud\`ene and Schapira extended Theorem \ref{T:Sigmund} to transitive flows on metric spaces with a standard Borel structure, provided that the flow admits a local product structure and satisfies the closing lemma; see \cite[Sections 2 and 3]{CS14}. In particular, their result applies to the geodesic flow on any complete connected negatively curved manifold of finite volume.

In this paper, we focus on the diagonal flow $(X_2,F)$, where $X_2:=\SL_2(\ZZ)\backslash\SL_2(\RR)$ and $F=\{a_t:t\in\RR\}$ is the connected diagonal subgroup of $\SL_2(\RR)$ acting on $X_2$ by right multiplication. This provides the homogeneous model for the geodesic flow on the unit tangent bundle of the modular surface.

Each $F$-periodic measure on $X_2$ has an associated real quadratic field $K$ and an order $\cR\subseteq K$. It can also be described by a rank-two lattice $L\subset K$, up to totally positive homothety. (See Section~\ref{SS:periodic-number} for details.) In this paper, we mainly use the data of the associated field $K$ and the discriminant of the suborder in $\cR$ generated by the expanding eigenvalue of a return matrix for the periodic orbit. Then it is natural to ask the following vague problem:
\begin{problem}\label{pro-main}
Is the set of periodic measures with prescribed arithmetic information still dense in the set of all invariant measures?
\end{problem}

One motivation for Problem \ref{pro-main} comes from the classical Duke's theorem presented by Einsiedler-Lindenstrauss-Michel-Venkatesh \cite{ELMV12}: if $\cG_d$ denotes the full packet of periodic geodesics with square-free discriminant $d>0$, then the normalized arc-length measure supported on $\cG_d$ converges weakly to the Haar measure on $X_2$ as $d\to +\infty$. Later, Aka and Einsiedler \cite{AE16} strengthened this result by proving the same equidistribution for every sequence of subcollections $\cI_d\subseteq\cG_d$ satisfying
$$
\frac{\ell(\mathcal G_d)}
{\ell(\mathcal I_d)\log d}\to 0,\quad \text{where }\ell(\cdot) \text{ denotes the total length of the collection}.
$$
However, these results only have the Haar measure as their target and need to take averages over collections of closed geodesics; neither of them provides an approximation theorem using measures supported on individual closed geodesics.

Concerning periodic measures supported on individual closed geodesics, McMullen \cite{M09} constructed, for every fixed real quadratic field, infinitely many primitive closed geodesics contained in a common compact subset of the modular surface. Furthermore, allowing the quadratic field to vary, Bourgain and Kontorovich \cite{BK17} proved that a fixed compact subset contains infinitely many fundamental closed geodesics, and obtained quantitative lower bounds for their abundance.

%Problem \ref{pro-main} is also motivated by recent works of McMullen \cite{M09} and Bourgain--Kontorovich \cite{BK17} on the distribution of periodic geodesics with prescribed arithmetic information.

In the fixed-field setting, Aka and Shapira \cite{AS18} proved that, for any real quadratic field $K$, there exists a sequence of individual periodic measures associated with $K$ that converges weakly to the Haar measure on $X_2$. More precisely, they established effective equidistribution of individual periodic orbits along certain Hecke rays.

Our goal in this paper is to address Problem \ref{pro-main} with certain arithmetic information on the discriminants of the approximating periodic measures. Let $\cS_f\subseteq\NN$ denote the set of square-free natural numbers, which is identified with the set of all real quadratic fields via the natural map $d\mapsto \QQ (\sqrt d)$. Our main theorem is:

\begin{theorem}\label{T:main1}
Let $A\subset \cS_f$ be a finite set of square-free numbers, and let $h: \NN\rightarrow \NN$ be a function with $\lim_{n\rightarrow +\infty} h(n)=+\infty$. Then there exists a subset $B\subset \cS_f\setminus A$ with
\begin{equation}\label{E:B-sparse}
	\limsup_{T\to +\infty} \frac{ |B \cap[0, T]|}{h(T)}=0,
\end{equation}
such that the set of all periodic measures with associated quadratic field from $B$ is dense in the set of all invariant measures.

Moreover, let $K\subset \SL_2(\ZZ)\backslash\SL_2(\RR)$ be a compact subset. Then there exists a compact subset $K'\supset K$, only depending on $K$, such that the set of all periodic measures that are supported in $K'$ and have associated quadratic field from $B$ is dense in the set of all invariant measures that are supported in $K$.
\end{theorem}

\begin{remark}
The requirement \eqref{E:B-sparse} means that the set $B=\{d_n:n\geq 1\}$ can be arbitrarily sparse. Moreover, we may additionally require that for all $n\geq 1$, $d_n\equiv 5\pmod{24}$; and that for all prime $p\geq 5$ with given signs $\epsilon_p\in \{\pm 1\}$ and all sufficiently large $n$, $\chi_p(d_n)=\epsilon_p$, where $\chi_p$ is the quadratic character on $\FF_p$. See Theorem \ref{T:main-arithmetic} and Remark \ref{R:strong} for this stronger conclusion.
\end{remark}

One may ask for a more restrictive condition than \eqref{E:B-sparse} on the admissible subset $B\subseteq\cS_f$. For example, given a finite subset $B\subseteq\cS_f$, it is interesting to investigate whether the set of all periodic measures with associated quadratic fields from $B$ is dense in the set of all invariant measures. % (see Section \ref{SS:rank-1-discuss} for more discussion). %However, using Hecke equidistribution, one can prove the following special case (see \cite[Theorem 4.8]{AS18}):
%\begin{theorem}
    %Let $d$ be a square-free number. Then the Haar measure on $X_2$ is a weak-${}^*$ limit of the set of all periodic measures with their associated discriminant all equal to $d$.  
%\end{theorem}
Furthermore, we raise the following problem: 

\begin{problem}\label{P:fixed-field-density}
Let $d\in\cS_f\setminus\{1\}$. Can every $F$-invariant measure on $X_2$ be approximated in the weak-${}^*$ topology by periodic measures whose associated quadratic fields are $\QQ(\sqrt d)$?
\end{problem}

\begin{remark}
Concerning Problem \ref{P:fixed-field-density}, a positive answer would give a fixed-field analogue of Sigmund's density theorem \cite[Theorem~1]{Sigmund72}. A negative answer would instead produce nontrivial arithmetic obstructions to approximation by periodic measures.
\end{remark}

Since periodic measures are dense in all invariant measures, Problem~\ref{P:fixed-field-density} amounts to asking the following: given $d\in\cS_f\setminus\{1\}$, can every periodic measure on $X_2$ be approximated by periodic measures associated with $\QQ(\sqrt d)$?

In particular, one may ask the following weaker question: given two distinct real quadratic fields $K_1$ and $K_2$, does there exist a periodic measure associated with $K_1$ that is a weak-${}^*$ limit of periodic measures associated with $K_2$?

The authors have no example or counterexample for any one of the above problems.

\subsection*{Organization of the paper} This paper is organized as follows. In Section \ref{S:pre}, we review some preliminary knowledge on the geodesic flow on $X_2$, including its arithmetic coding and its number-theoretic correspondence. In Section \ref{S:cpt-shadow}, we prove an orbit shadowing result for the geodesic flow in compact subsets of $X_2$, which is the first technical contribution of this paper. In Section \ref{S:periodic-arithmetic}, we establish a general proposition that allows us to add certain arithmetic corrections on the periodic approximation of invariant measures.  

The main novelty of this paper lies in converting the problem of periodic approximation with prescribed arithmetic information into the corresponding problem of modifying the tails of periodic words in the full-shift system. This conversion heavily relies on the precise translation from the language of invariant measures to the world of continued fractions. Throughout Sections \ref{S:pre}, \ref{S:cpt-shadow}, and \ref{S:periodic-arithmetic}, we repeatedly highlight the significance of this translation and emphasize its central role in our approach.

In Section \ref{S:arithmetic-density}, we prove a stronger result from which our main theorem (Theorem \ref{T:main1}) follows. In Section \ref{S:discuss}, we discuss the higher-rank analogues of Theorem \ref{T:main1} and Problem \ref{P:fixed-field-density}, including their connection with the four exponentials conjecture and the six exponentials theorem.

\subsection*{Acknowledgment}  
The authors are grateful to Antoine Marnat for helpful discussions, through which they learned about the four exponentials conjecture and the six exponentials theorem in Section \ref{S:discuss}.

This work is supported by National Key R\&D Program of China No. 2022YFA1007800, NSFC No. 12071474 and Zhejiang Provincial Innovative Resource Allocation Project  2022XHSJJ010.

\subsection*{Declaration of AI use} Indeed, the main results of this paper were motivated by the higher-rank case discussed in Section \ref{S:discuss}, which was known to the authors since 2016. 
The current project was initiated in early 2025. The authors independently developed all mathematical ideas in the proofs of their results, and had given two talks in Xiamen University and Soochow University, at the end of 2025 and beginning of 2026 respectively. Later, with the assistance of AI tools, the authors refined certain estimates in Section \ref{S:cpt-shadow} and simplified several arguments in Section \ref{S:periodic-arithmetic}. The authors have carefully verified all calculations and proofs and take full responsibility for the entire manuscript.

\section{Preliminaries}\label{S:pre}

\subsection{Basics on the Gauss map and the geodesic flow}
In this subsection, we review some basic knowledge on the invertible Gauss map and its suspension into the geodesic flow on $X_2=\SL_2(\ZZ)\backslash\SL_2(\RR)$. More detailed discussion can be found in \cite{AS13,GK01,Series85}.

\subsubsection{An invertible extension of the usual Gauss map}

Let $\Sigma=\NN^{\ZZ}$ be the full shift space over the countable alphabet $\NN=\{1,2,\cdots\}$, and let $\sigma$ denote the left shift on $\Sigma$:
$$
\ba=(a_j)_{j\in\ZZ}\mapsto\sigma(\ba)=(a_{j+1})_{j\in\ZZ}.
$$
For $\ba\in\Sigma$, we associate it with two continued fractions $$
x(\ba)=[0;a_1,a_2,\cdots],\quad \text{and}\quad y(\ba)=[0;a_0,a_{-1},\cdots].
$$
Then $x(\ba),y(\ba)\in (0,1)\setminus\QQ$. We also introduce $$
z(\ba)=\frac{y(\ba)}{1+x(\ba)y(\ba)}.
$$
Then $0<z(\ba)<\frac{1}{1+x(\ba)}$. Moreover, it is clear that $$
x(\sigma(\ba))=\frac{1}{x(\ba)}-\left\lfloor \frac{1}{x(\ba)}\right\rfloor,\quad \text{and}\quad
z(\sigma(\ba))=x(\ba)\cdot (1-x(\ba)z(\ba)).
$$

Let
\[
\begin{aligned}
D&:=\left\{(x,z)\in\RR^2:0<x<1,\ 0<z<\frac{1}{1+x}\right\},\\
D^*&:=\left\{(x,z)\in D:x\notin\QQ,\ \frac{z}{1-xz}\notin\QQ\right\},
\end{aligned}
\]
and consider the map
$$
T\colon D^*\to D^*,\quad (x,z)\mapsto \left(\frac{1}{x}-\left\lfloor \frac{1}{x}\right\rfloor, x(1-xz)\right).
$$
It follows that the coding map $$
\Sigma\to D^*,\quad \ba\mapsto (x(\ba),z(\ba))
$$
is a bijection which conjugates $\sigma$ with $T$. So the dynamical system $(\Sigma,\sigma)\cong (D^*,T)$ is an invertible extension of the usual Gauss map $$
(0,1)\setminus\QQ\to (0,1)\setminus\QQ,\quad x\mapsto \frac{1}{x}-\left\lfloor \frac{1}{x}\right\rfloor.
$$

\subsubsection{Two chambers of the homogeneous space $X_2$}

Let $G=\SL_2(\RR)$, $\Gamma=\SL_2(\ZZ)$. We explicitly construct two disjoint chambers of $X_2$ whose union almost equals $X_2$.

We write \begin{align*}
    \Omega_-:=\left\{
\begin{pmatrix}
    g_{11} &g_{12}\\ g_{21} &g_{22}
\end{pmatrix}\in G: 
0<-g_{21}<1\leq g_{11},\;0<g_{12}<g_{22}
\right\};\\
\Omega_+:=\left\{
\begin{pmatrix}
    g_{11} &g_{12}\\ g_{21} &g_{22}
\end{pmatrix}\in G: 
0<g_{11}<1\leq -g_{21},\;0<g_{22}<g_{12}
\right\}.
\end{align*}
The diagonal flow is given by right multiplication by $$
a_t:=\begin{pmatrix}
    e^{t/2} &0\\ 0 &e^{-t/2}
\end{pmatrix}\,(t\in\RR).
$$
Let $\{e_1,e_2\}$ denote the standard basis of $\RR^2$, and consider the following full subset of $X_2$: \begin{equation}\label{E:X2-0}
X_2^0:=\{\Gamma g\in X_2:\ZZ^2g\cap (\RR e_1\cup\RR e_2)=\{0\}\}.
\end{equation}
It is readily checked that each point in $X_2^0$ can be uniquely written as $\Gamma g$ where $g\in \Omega_{\pm}$. In particular, $$X_2^0=(\Gamma\Omega_-\cap X_2^0)\sqcup(\Gamma\Omega_+\cap X_2^0).$$

Moreover, for any $(x,z)\in (0,1)^2$, we write $$
g_-(x,z)=\begin{pmatrix}
    1 &z\\ -x &1-xz
\end{pmatrix},\quad 
g_+(x,z)=\begin{pmatrix}
    x &1-xz\\ -1 &z
\end{pmatrix}.
$$
It is readily checked that \begin{align*}
    \Omega_-&=\left\{g_-(x,z)\cdot a_t:(x,z)\in D,0\leq t<-2\log{x}\right\};\\
    \Omega_+&=\left\{g_+(x,z)\cdot a_t:(x,z)\in D,0\leq t<-2\log{x}\right\}.
\end{align*}
We conclude that each point in $X_2^0$ can be uniquely written as $\Gamma g_{\pm}(x,z)\cdot a_t$ for some $(x,z)\in D^*$ and some $0\leq t<-2\log{x}$.
The incoming boundaries of the chambers $\Gamma\Omega_{\pm}$ are \begin{equation}
    \cC_{\pm}:=\{\Gamma g_{\pm}(x,z):(x,z)\in D\}.
\end{equation}
They will be regarded as the cross-sections in $X_2$ that are transversal to $\{a_t:t\in\RR\}$-actions.

In particular, for $x\in (0,1)$ we write $\tau(x)=-2\log{x}$. Then for any $(x,z)\in D^*$, $$
\Gamma g_{\epsilon}(x,z)\cdot a_{\tau(x)}=\Gamma g_{-\epsilon}(T(x,z)).
$$

\subsubsection{A suspension of the invertible Gauss map into the geodesic flow}
Let $\widehat{\Sigma}=\Sigma\times \{\pm 1\}$ be the two-folded space of $\Sigma$, and let $\widehat{\sigma}$ denote the two-folded left shift on $\widehat{\Sigma}$: $$
(\ba,\epsilon)\mapsto (\sigma(\ba),-\epsilon).
$$
We also write $$
\cC_{\pm}^*:=\{\Gamma g_{\pm}(x,z):(x,z)\in D^*\},
$$
and consider the return map $$
\widehat{T}\colon\cC_{\pm}^*\to \cC_{\pm}^*,\quad \Gamma g_{\epsilon}(x,z)\mapsto \Gamma g_{\epsilon}(x,z)\cdot a_{\tau(x)}.
$$
It follows that the coding map $$
\widehat{\Sigma}\to \cC_{\pm}^*, \quad (\ba,\epsilon)\mapsto \Gamma g_{\epsilon}(x(\ba),z(\ba))
$$
is a bijection which conjugates $\widehat{\sigma}$ with $\widehat{T}$. This gives a suspension of the two-folded invertible Gauss map into the geodesic flow on $X_2$ with cross-sections $\cC^*_{\pm}$: \begin{equation}\label{E:susmap}
    \Psi\colon((\ba,\epsilon),t)\mapsto \Gamma g_{\epsilon}(x(\ba),z(\ba))\cdot a_t,\quad \text{where}\quad 0\leq t<\tau(x(\ba)).
\end{equation}
In particular, any even-length periodic digit word $\ba\in\Sigma$ with a choice of chamber $\epsilon\in\{\pm1\}$ gives a closed geodesic on $X_2$.

\subsubsection{An interpretation of cross-sections in the upper-half plane} We give a geometric interpretation of the cross-sections constructed above in terms of geodesic flows on the upper-half plane.

Let $$\HH:=\{z\in\CC:\mathrm{Im}(z)>0\}$$ be equipped with the standard hyperbolic metric. It is well known that the action of $G$ on $\HH$ induced by fractional linear transformations $$
G\times \HH\to\HH,\quad 
(g,z)\mapsto g.z:=\frac{g_{11}z+g_{12}}{g_{21}z+g_{22}}
$$
preserves the metric. This action furthermore extends to an action of $G$ on the unit tangent bundle $\rT^1\HH$: \begin{equation}\label{E:Gaction}
G\times \rT^1\HH\to\rT^1\HH,\quad 
(g,(z,\zeta))\mapsto (g.z,Dg(\zeta)).
\end{equation}
We see that the action \eqref{E:Gaction} is transitive and has kernel $\{\pm 1_G\}$, so it gives the explicit isomorphism \begin{equation*}
G/\{\pm 1_G\}\xrightarrow{\cong} \rT^1\HH,\quad g\mapsto (g.i,Dg(\uparrow)).
\end{equation*} Moreover, under this isomorphism, the right multiplication by $\{a_t:t\in\RR\}$ on $G$ corresponds to the geodesic flow on $\rT^1\HH$.

Let us explicitly calculate the geodesic on $\HH$ through $g.i$ where $g\in\Omega_{\pm}$. First, the expression $$
\begin{pmatrix}
    g_{11} &g_{12}\\ g_{21} &g_{22}
\end{pmatrix}\cdot 
\begin{pmatrix}
    e^{t/2} &0\\ 0 &e^{-t/2} 
\end{pmatrix}.i
$$
converges to $r_2(g)=\frac{g_{12}}{g_{22}}\in\RR$ as $t\to -\infty$, and converges to $r_1(g)=\frac{g_{11}}{g_{21}}\in\RR$ as $t\to +\infty$. So the geodesic through $g.i$ is a half-circle in $\HH$ connecting $r_2(g)$ to $r_1(g)$. Moreover, it is readily checked that \begin{align*}
    \Omega_-&=\left\{g_-\left(-\frac{1}{r_1},-\frac{r_1r_2}{r_2-r_1}\right)\cdot a_t:r_1<-1,0<r_2<1,0\leq t<-2\log(-r_1^{-1})\right\};\\
    \Omega_+&=\left\{g_+\left(-r_1,\frac{1}{r_2-r_1}\right)\cdot a_t:-1<r_1<0,r_2>1,0\leq t<-2\log(-r_1)\right\}.
\end{align*}

\subsection{Basics on the $F$-periodic orbits in $X_2$}\label{SS:periodic-number}
In this subsection, we review some basic knowledge on the construction and description of $F$-periodic orbits in $X_2$. One may refer to \cite[Sections~1--2]{M09} and \cite[Sections~2.2--2.3]{ELMV12} for more backgrounds and details.

\subsubsection{Construction of $F$-periodic measures on $X_2$}\label{SSS:construct-periodic} Recall that $$\Psi\colon((\ba,\epsilon),t)\mapsto \Gamma g_{\epsilon}(x(\ba),z(\ba))\cdot a_t,\quad \text{where}\quad 0\leq t<\tau(x(\ba))$$
is a suspension of the two-folded invertible Gauss map into the geodesic flow on $X_2$ with cross-sections $\cC^*_{\pm}$. Now we fix $\epsilon=+1$ and consider one branch of the suspension: \begin{equation}\label{E:+susmap}
    \Psi_+\colon \Sigma\times\RR\to X_2,\qquad (\ba,t)\mapsto \Gamma g_{+}(x(\ba),z(\ba))\cdot a_t.
\end{equation}
Recall that for $x\in (0,1)$, $\tau(x)=-2\log{x}$. For $\ba\in\Sigma$, we define its $+$ return time to be $$
R(\ba):=\tau(x(\ba))+\tau(x(\sigma(\ba))).
$$
It is easy to see that $$
R(\ba)=-2\log{\frac{x(\sigma(\ba))}{a_1+x(\sigma(\ba))}}>2\log{2}.
$$
Moreover, we have $\Psi_+(\ba,R(\ba)+t)=\Psi_+(\sigma^2(\ba),t)$ for all $(\ba,t)\in\Sigma\times\RR$. This means that $\Psi_+$ is the suspension of the return map $\widehat{\sigma}^2$ on the $+$ chamber $\Sigma\times\{+1\}$ of $\widehat{\Sigma}$. 

For any bounded Borel function $f:X_2\to\RR$ and any $\ba\in \Sigma$, write \begin{equation}
    J_f(\ba):=\int_0^{R(\ba)}f(\Psi_+(\ba,t))\mathrm{d}t.
\end{equation}
This is the integral of $f$ over the first-return orbit segment $\{\Psi_+(\ba,t):0\leq t\leq R(\ba)\}$. 

Now we use the suspension \eqref{E:+susmap} to construct $F$-periodic orbits on $X_2$. Given a nonempty word $W$ of even length $2k$ over the alphabet $\NN$, we denote its bi-infinite repetition by $\overline{W}\in\Sigma$; that is, for $W=a_1a_2\cdots a_{2k-1}a_{2k}$, we write $$
\overline{W}=(\cdots;a_1,a_2,\cdots,a_{2k-1},a_{2k};a_1,a_2,\cdots,a_{2k-1},a_{2k};\cdots).
$$
For any bounded Borel function $f:X_2\to\RR$, define \begin{equation}
    I_f(W):=\sum_{j=0}^{k-1}J_f(\sigma^{2j}(\overline{W})).
\end{equation}
This is the integral of $f$ over the parametrized closed orbit $$
\cO_W:=\left\{\Psi_+(\overline{W},t):
0\leq t\leq \sum_{j=0}^{k-1}R(\sigma^{2j}(\overline{W}))
\right\}.
$$
In particular, we see that the full period of this parametrized closed orbit is $$
T(W):=I_1(W)=\sum_{j=0}^{k-1}R(\sigma^{2j}(\overline{W})).
$$
It follows that $T(W)\geq k\cdot 2\log{2}=|W|\cdot\log{2}$.

Finally, we construct $F$-periodic measures on $X_2$ as follows. Given a nonempty word $W$ of even length $2k$ over the alphabet $\NN$, we define a measure $\mu_W$ on $X_2$ by \begin{equation}
\mu_W(f):=\frac{I_f(W)}{T(W)},\qquad \text{for all }f\in C_b(X_2).
\end{equation}
This is the normalized uniform measure on the closed orbit $\cO_W$.

\subsubsection{Correspondence between closed orbits in $X_2$ and quadratic fields}\label{SSS:correspondence}

We recall the arithmetic information carried by a periodic orbit of the diagonal flow on $X_2$. Let $\cO=\Gamma gF$ be an periodic orbit with a period $T>0$, and write $F=\{a_t:t\in\RR\}$. Then $\gamma:=g a_Tg^{-1}\in\Gamma.$
Write
\[
 \Delta(\gamma):=(\operatorname{tr}\gamma)^2-4>0,\qquad
 K_{\mathcal O}:=\mathbb Q\bigl(\sqrt{\Delta(\gamma)}\bigr).
\]
Since $\Delta(\gamma)$ is not a square, $K_{\cO}$ is a real quadratic field. This field is independent of the choice
of $g$ and also independent of the positive period $T$. Moreover, there is a field isomorphism $$
\Phi_{\gamma}:\QQ[\gamma]\to K_{\cO},\qquad aI_{2}+b\gamma\mapsto a+b\lambda \quad (a,b\in\QQ),\quad \text{where }\lambda:=e^{T/2}>1.
$$

The orbit carries finer arithmetic information than the field alone. Write $\Lambda=\mathbb Z^2g$.
%, and
%\[
 %\End_{\mathrm{diag}}(\Lambda):=\left\{D\in\RR^{2\times 2}\text{ diagonal}:\Lambda\cdot D\subseteq\Lambda\right\}.
%\]
%Then there is a ring isomorphism $$\End_{\mathrm{diag}}(\Lambda)\to \QQ[\gamma]\cap \ZZ^{2\times 2},\qquad D\mapsto gDg^{-1}.$$
Then $$\mathcal R_{\mathcal O}:=\Phi_{\gamma}(\QQ[\gamma]\cap \ZZ^{2\times 2})$$ is an order in $K_{\mathcal O}$. This order is an invariant of the orbit and it satisfies
$\mathbb Z[\lambda]\subseteq \mathcal R_{\mathcal O}\subseteq \mathcal O_{K_{\mathcal O}},$
where $\mathcal O_{K_{\mathcal O}}$ is the maximal order of
$K_{\mathcal O}$. In particular,
\[
 \operatorname{disc}\bigl(\mathbb Z[\lambda]\bigr)=\Delta(\gamma).
\]
%it need not be the discriminant of $\mathcal R_{\mathcal O}$ or of $K_{\mathcal O}$.

Conversely, let $K$ be a real quadratic field with %(fixed ordered) --it doesn't matter
real embeddings $\sigma_1,\sigma_2$ and $\mathcal R\subset K$ an order.
We choose a rank-two lattice $I\subset K$ with $\{\alpha\in K:\alpha\cdot I\subseteq I\}=\mathcal R$. Consider the lattice
\[
 \sigma(I):=\{(\sigma_1(x),\sigma_2(x)):x\in I\}\subset\RR^2,
\]
and normalize it to have covolume one:
\[
 \Lambda_I:=\operatorname{covol}(\sigma(I))^{-1/2}\cdot\sigma(I).
\]
Then $\Lambda_I$ represents a point of $X_2$ whose $F$-orbit is periodic. The associated quadratic field and multiplier order of this orbit are $K$ and $\mathcal R$, respectively.

Two such lattices $I,J\subset K$ determine the same $F$-orbit if and only if $J=\alpha I$ for some totally positive $\alpha\in K^\times$. Consequently, for a fixed quadratic field $K$ and an order $\mathcal R$, the corresponding periodic $F$-orbits are parametrized by the narrow ideal classes of $\mathcal R$. %See \cite[Sections~2.2--2.3]{ELMV12} for the underlying correspondence between quadratic orders, optimal embeddings, and ideal classes.

\subsubsection{Connection between continued-fractions and discriminants}\label{SSS:discriminant}
For $n\in\NN$, write $$\sff(n)=d \qquad \text{for } n=ds^2 \text{ with } d \text{ square-free and }s\in\NN.$$
For $a\in\NN$, define
\[
 C(a):=\begin{pmatrix}
 a&1\\ 1&0
 \end{pmatrix}.
\]
For a finite word $W=a_1\cdots a_\ell$ over the alphabet $\NN$, put
\[
 C(W):=C(a_1)\cdots C(a_\ell),\qquad C(\varnothing):=I.
\]
Since $\det C(a)=-1$, we have $\det C(W)=(-1)^\ell.$ In particular, $C(W)\in\SL_2(\ZZ)$ whenever $W$ has even length.

Let $W$ be a nonempty word of even length. We define its matrix
discriminant by
\begin{equation}\label{E:word-discriminant}
 \Delta(W):=(\operatorname{tr}C(W))^2-4
           =d(W)s(W)^2,
 \qquad d(W):=\sff(\Delta(W)).
\end{equation}
The next lemma identifies the quadratic field associated with
the orbit $\cO_W$ constructed in
Section~\ref{SSS:construct-periodic}.

\begin{lemma}\label{L:word-field}
Let $W$ be a nonempty word of even length. Then
$\Delta(W)$ is a positive nonsquare integer, and the quadratic
field associated with $\cO_W$ is
\[
 K_{\cO_W}=\QQ\bigl(\sqrt{\Delta(W)}\bigr)
          =\QQ\bigl(\sqrt{d(W)}\bigr).
\]
In particular, $d(W)>1$.
\end{lemma}
\begin{proof}
Write $W=a_1a_2\cdots a_{2k-1}a_{2k}$, and put
\[
 x_i:=x(\sigma^{i-1}(\overline W))\quad(1\leq i\leq 2k),\qquad x_{2k+1}:=x_1.
\]
Recall that the parametrized closed orbit $$
\cO_W=\left\{\Psi_+(\overline{W},t):
0\leq t\leq \sum_{j=0}^{k-1}R(\sigma^{2j}(\overline{W}))
\right\}.
$$
has its full period $$
T(W)=\sum_{j=0}^{k-1}R(\sigma^{2j}(\overline{W}))=-2\sum_{j=0}^{k-1}(\log{x_{2j+1}}+\log{x_{2j+2}})
=-2\sum_{i=1}^{2k}\log{x_i}.
$$
This means that the closed orbit $\cO_W$ can be written as $\Gamma g_WF$, where $$
g_W=g_+(x(\overline{W}),z(\overline{W}))\in G,\quad\text{and}\quad
\gamma_W:=g_Wa_{T(W)}g_W^{-1}\in\Gamma.
$$
We claim that: \begin{equation}\label{E:same-tr}
    \tr{\gamma_W}=\tr{C(W)}.
\end{equation}

In fact, write $\lambda_W=\prod_{i=1}^{2k}x_i^{-1}>1$. It follows that $$\tr{\gamma_W}=\tr{a_{T(W)}}=e^{T(W)/2}+e^{-T(W)/2}=\lambda_W+\lambda_W^{-1}.$$
On the other hand, for $1\leq i\leq 2k$, since $x_i=\frac{1}{a_i+x_{i+1}}$, we have
\[
 C(a_i)\begin{pmatrix}1\\x_{i+1}\end{pmatrix}=
 x_i^{-1}\begin{pmatrix}1\\x_i\end{pmatrix}.
\]
Using these identities for $i=1,\cdots,2k$ with $x_{2k+1}=x_1$, we obtain
\[
 C(W)\begin{pmatrix}1\\x_1\end{pmatrix}
 =\lambda_W\begin{pmatrix}1\\x_1\end{pmatrix}.
\]
Since $\det C(W)=1$, the eigenvalues of $C(W)$ are $\lambda_W$ and $\lambda_W^{-1}$. This gives $\tr{C(W)}=\lambda_W+\lambda_W^{-1}$ and hence proves \eqref{E:same-tr}.

Therefore, by the correspondence in Section \ref{SSS:correspondence},
\[
 K_{\cO_W}=\QQ\bigl(\sqrt{(\tr\gamma_W)^2-4}\bigr)=\QQ\bigl(\sqrt{(\tr{C(W)})^2-4}\bigr)
 =\QQ\bigl(\sqrt{\Delta(W)}\bigr).
\]
Note that $\tr C(W)\geq 3$ since it is an integer greater than $2$. Thus $\Delta(W)=t^2-4$ is positive and is not a square. This implies that $d(W)>1$.
\end{proof}

We emphasize that $\Delta(W)$ is the discriminant of the order $\ZZ[\lambda_W]$. It is not necessarily the discriminant of the multiplier order of $\cO_W$, nor the discriminant of its quadratic field. If $W$ is replaced by a repetition $W^m$, the underlying closed orbit remains the same, while $\Delta(W^m)$ may change. Nevertheless, Lemma~\ref{L:word-field} shows that the square-free part $d(W)$, and hence the associated quadratic field, remains unchanged.

For later use, let $p$ be an odd prime. The quadratic character is defined by
\[
\chi_p:\FF_p\to \{0,\pm1\},\qquad x\mapsto \begin{cases}
  0,  & x=0,\\
  1,  & x\in(\FF_p^{\times})^2,\\
  -1, & x\notin(\FF_p^{\times})^2\text{ and }x\neq 0.
 \end{cases}
\]
For $n\in\ZZ$, we also write $\chi_p(n)=\chi_p(n\mod{p})$. It follows that for a nonempty word $W$ of even length, $$
\chi_p(d(W))=
\begin{cases}
0,&p\text{ ramifies}\\
1,&p\text{ splits}\\
-1,&p\text{ is inert}\\
\end{cases}
\quad\text{in the field } K_{\cO_W}=\QQ(\sqrt{d(W)}).
$$

We write $v_p(n)$ for the exponent of $p$ in a nonzero integer $n$. In particular, if $v_p(\Delta(W))=1$, then $p$ ramifies in $K_{\cO_W}$; if $v_p(\Delta(W))=0$, then $p$ either splits or is inert in $K_{\cO_W}$.

\section{A compact shadowing estimate}\label{S:cpt-shadow}

In this section, we prove an orbit shadowing result in compact subsets of $X_2$ provided that the intersections of these orbits with the cross-sections are close to each other.

\subsection{Exponentially close data of continued fractions} Write $$
\theta:=\frac{\sqrt5-1}{2}\in (0,1).
$$
We begin with a simple lemma concerning the distance of two continued fractions with the same first $n$ partial quotients.

\begin{lemma}\label{L:samepartial}
    There exists an absolute constant $C_0>0$ such that two continued fractions with the same first $n$ partial quotients differ by at most $C_0\theta^{2n}$.
\end{lemma}
\begin{proof}
    Let $\frac{p_n}{q_n}$ be the common $n$-th convergent of the two continued fractions $\alpha,\beta$. Then we have $$
    |\alpha-\beta|\leq \left|\alpha-\frac{p_n}{q_n}\right|+\left|\beta-\frac{p_n}{q_n}\right|\leq \frac{2}{q_n^2}.
    $$
    Note that $q_n\geq F_n$ where $F_n$ is the $n$-th Fibonacci number. Since $F_n\approx \theta^{-n}$, the conclusion follows.
\end{proof}

\begin{corollary}\label{C:samepartial}
    There exists an absolute constant $C_1>0$ such that the following holds. For $\ba,\ba'\in\Sigma$ whose digit numbers agree on $[-n,n]$, where $n\ge1$, one has \begin{align*}
        |x(\ba)-x(\ba')|+|y(\ba)-y(\ba')|+|z(\ba)-z(\ba')|&\leq C_1\theta^{2n};\\
        \left|\frac{1}{x(\ba)}-\frac{1}{x(\ba')}\right|+\left|\frac{1}{y(\ba)}-\frac{1}{y(\ba')}\right|&\leq C_1\theta^{2n};\\
        |\tau(x(\ba))-\tau(x(\ba'))|&\leq C_1\theta^{2n}.
    \end{align*}
\end{corollary}
\begin{proof}
    Since the map $$
    (0,1)^2\to (0,1),\quad (x,y)\mapsto z=\frac{y}{1+xy}
    $$
    has bounded derivatives, the first estimate follows from Lemma \ref{L:samepartial}. Note that $$
    \frac{1}{x(\ba)}=a_1+x(\sigma(\ba))\quad \text{and}\quad \frac{1}{y(\ba)}=a_0+y(\sigma^{-1}(\ba)).
    $$
    The second estimate also follows from Lemma \ref{L:samepartial}. Finally, since $$
    \tau(x(\ba))=2\log{\frac{1}{x(\ba)}}\quad \text{where}\quad \frac{1}{x(\ba)}>1,
    $$
    and the function $\log(\cdot)$ is $1$-Lipschitz on $(1,+\infty)$, the third estimate follows from the second one.
\end{proof}

\subsection{Distance between geodesic flows on bounded intervals} In view of the suspension map \eqref{E:susmap}, we may control the distance between two geodesic flows using their intersections with the cross-sections $\cC_{\pm}$ and their travel time after the intersections. 

We first fix some notations on $X_2=\Gamma\backslash G=\SL_2(\ZZ)\backslash\SL_2(\RR)$. The Frobenius norm $\|\cdot\|$ on the Lie algebra of $G$ determines a left-invariant metric $d_G$ on $G$, and hence induces a Riemannian metric $d_{X_2}$ on $X_2$.

\begin{lemma}\label{L:Lipdist}
    For any $T>0$, any $(x,z),(x',z')\in (0,1)^2$, any $\epsilon\in\{\pm 1\}$, and any $t,t'\in [-T,T]$, one has \begin{equation}
        d_{X_2}(\Gamma g_{\epsilon}(x,z)a_t,\Gamma g_{\epsilon}(x',z')a_{t'})\leq 2e^T(|x-x'|+|z-z'|)+|t-t'|.
    \end{equation}
\end{lemma}
\begin{proof}
Let $\epsilon\in \{\pm 1\}$. For $p=(x,z;t)\in (0,1)^2\times [-T,T]$, we write $$
g_{\epsilon}(p)=g_{\epsilon}(x,z)\cdot a_t.
$$
It is clear that $$
\varphi:[0,1]\to G,\quad s\mapsto g_{\epsilon}(p_s)\quad \text{where}\quad p_s:=(1-s)p+sp'
$$
is a continuous curve in $G$ connecting $g_{\epsilon}(p)$ to $g_{\epsilon}(p')$. It follows that \begin{equation*}
    d_{X_2}(\Gamma g_{\epsilon}(x,z)a_t,\Gamma g_{\epsilon}(x',z')a_{t'})\leq d_G(g_{\epsilon}(p),g_{\epsilon}(p'))\leq \int_0^1\|\varphi'(s)\|_{\varphi(s)}\mathrm{d}s.
\end{equation*}
It suffices to calculate the norm $\|\varphi'(s)\|_{\varphi(s)}$.

Note that the chain rule gives $$
\varphi^{-1}(s)\varphi'(s)=g_{\epsilon}(p_s)^{-1}\cdot 
\left(\frac{\partial g_{\epsilon}}{\partial x}(p_s)\cdot (x'-x)+\frac{\partial g_{\epsilon}}{\partial z}(p_s)\cdot (z'-z)+\frac{\partial g_{\epsilon}}{\partial t}(p_s)\cdot (t'-t)\right).
$$
Write $p_s=(x_s,z_s;t_s)$. We calculate that
\begin{equation*}
    \begin{aligned}
    g_{\epsilon}(p_s)^{-1}\cdot\frac{\partial g_{\epsilon}}{\partial x}(p_s)
    &=\begin{pmatrix}
    z_s &-\epsilon e^{-t_s}z_s^2\\ \epsilon e^{t_s} &-z_s
    \end{pmatrix};\\
    g_{\epsilon}(p_s)^{-1}\cdot\frac{\partial g_{\epsilon}}{\partial z}(p_s)
    &=\begin{pmatrix}
    0 &-\epsilon e^{-t_s}\\ 0 &0
    \end{pmatrix};\\
    g_{\epsilon}(p_s)^{-1}\cdot\frac{\partial g_{\epsilon}}{\partial t}(p_s)
    &=\begin{pmatrix}
    \frac{1}{2} &0\\ 0 &-\frac{1}{2}
    \end{pmatrix}.
\end{aligned}
\end{equation*}
So we have $$
\|\varphi'(s)\|_{\varphi(s)}=\|\varphi(s)^{-1}\varphi'(s)\|_{1_G}
\leq 2e^{T}(|x-x'|+|z-z'|)+|t-t'|.
$$
Integrating this over $s\in [0,1]$ gives the desired estimate.
\end{proof}

To control the travel time of the geodesics after its intersection with the cross-sections, we introduce a cusp-height function on $X_2$ and give a quantitative estimate of heights of the suspension map \eqref{E:susmap}. Define \begin{equation}\label{E:ht}
    \het: X_2\to [0,+\infty),\quad \Gamma g\mapsto \sup_{\gamma\in\Gamma}\mathrm{Im}(\gamma g.i).
\end{equation}

\begin{lemma}\label{L:compactdist}
For any $(\ba,\epsilon)\in\widehat{\Sigma}$ and any $0\leq t<\tau(x(\ba))$, one has \begin{equation}
    \het(\Psi((\ba,\epsilon),t))\geq \frac{1}{2}\exp\left(\min\{t,\tau(x(\ba))-t\}\right).
\end{equation}
In particular, let $K\subseteq X_2$ be a compact subset and write 
\begin{equation}\label{def-tk}
T_K:=\max\left\{0,\;\log\left(2\max_{q\in K}\het(q)\right)\right\} \text{ for } \het:X_2\to [0,+\infty)\text{ is defined in \eqref{E:ht}}.
\end{equation}
Then for any $0\leq t<\tau(x(\ba))$, \begin{equation}
    \Psi((\ba,\epsilon),t)\in K\Longrightarrow \min\{t,\tau(x(\ba))-t\}\leq T_K.
\end{equation}
\end{lemma}
\begin{proof}
For $g=\begin{pmatrix}
    g_{11} &g_{12}\\ g_{21} &g_{22}
\end{pmatrix}$, one has $\mathrm{Im}(g.i)=(g_{21}^2+g_{22}^2)^{-1}.$ Let $\lambda_1(\ZZ^2g)$ denote the length of the shortest nonzero vector in the row lattice $\ZZ^2g$. Then we have $$
\het(\Gamma g)=\lambda_1(\ZZ^2g)^{-2}.
$$
In particular, the cusp-height function is finite and continuous, so for any compact subset $K\subseteq X_2$, one has $\max_{q\in K}\het(q)<+\infty$.

For $\epsilon=-1$, we directly calculate that $$
\mathrm{Im}(g_-(x,z)\cdot a_t.i)=\frac{1}{e^tx^2+e^{-t}(1-xz)^2};
$$
for $\epsilon=+1$, we have $$
\mathrm{Im}\left(\begin{pmatrix}
    0 &1\\-1 &0
\end{pmatrix}\cdot g_+(x,z)\cdot a_t.i\right)=\frac{1}{e^tx^2+e^{-t}(1-xz)^2}.
$$
So in any case, \begin{equation}
\het(\Psi((\ba,\epsilon),t))\geq\frac{1}{e^tx(\ba)^2+e^{-t}(1-x(\ba)z(\ba))^2}.
\end{equation}
Note that $$
e^tx(\ba)^2=e^{t-\tau(x(\ba))},\quad \text{and}\quad e^{-t}(1-x(\ba)z(\ba))^2\leq e^{-t}.
$$
We conclude that $$
\het(\Psi((\ba,\epsilon),t))\geq\frac{1}{2}\exp\left(\min\{\tau(x(\ba))-t,t\}\right).
$$
The last statement clearly follows from this.
\end{proof}
    
Finally, we summarize all discussions in the section into a compact shadowing result.

\begin{proposition}\label{P:cpt-shadow}
    There exists an absolute constant $n_0>0$ such that the following holds. For any compact subset $K\subseteq X_2$, there exists a constant $C_K>0$ such that for $\ba,\ba'\in\Sigma$ whose digit numbers agree on $[-n,n]$ where $n\geq n_0$, if $$
    0\leq t<\min\{\tau(x(\ba)),\tau(x(\ba'))\}\quad \text{and}\quad \Psi((\ba,\epsilon),t)\in K,
    $$
    then $$
    d_{X_2}(\Psi((\ba,\epsilon),t),\Psi((\ba',\epsilon),t))\leq C_K\theta^{2n}.
    $$
\end{proposition}
\begin{proof}
    By Corollary \ref{C:samepartial}, there exists an absolute constant $C_1>0$ such that $$
    |x(\ba)-x(\ba')|+|z(\ba)-z(\ba')|\leq C_1\theta^{2n},\quad \text{and}\quad |\tau(x(\ba))-\tau(x(\ba'))|\leq C_1\theta^{2n}.
    $$
    Let $T_K\ge 0$ be the constant in Lemma \ref{L:compactdist}. Then our assumption on $t$ implies that $$
    \min\{t,\tau(x(\ba))-t\}\leq T_K.
    $$

    If $t\leq T_K$, then it follows from Lemma \ref{L:Lipdist} that $$
    d_{X_2}(\Psi((\ba,\epsilon),t),\Psi((\ba',\epsilon),t))\leq 2e^{T_K}(|x(\ba)-x(\ba')|+|z(\ba)-z(\ba')|)\leq 2e^{T_K}C_1\theta^{2n}.
    $$
    
    If $\tau(x(\ba))-t\leq T_K$, then we write $$
    \widetilde{t}:=t-\tau(x(\ba))\quad \text{and}\quad \widetilde{t}':=t-\tau(x(\ba')).
    $$
    It is clear that $\widetilde{t}\in [-T_K,0]$ and $|\widetilde{t}-\widetilde{t}'|\leq C_1\theta^{2n}$. In particular, when $n\geq n_0$ is large enough, one has $\widetilde{t}'\in [-T_K-1,0]$. We write $\Psi((\ba,\epsilon),t)=\Gamma g_{\epsilon}(x(\ba),z(\ba))\cdot a_t$ for any $(\ba,\epsilon)\in\widehat{\Sigma}$ and any $t\in\RR$, and note that $$
    \Psi((\ba,\epsilon),t)=\Psi((\sigma(\ba),-\epsilon),\widetilde{t})\quad\text{and}\quad \Psi((\ba',\epsilon),t)=\Psi((\sigma(\ba'),-\epsilon),\widetilde{t}').
    $$
    It follows from Lemma \ref{L:Lipdist} again that \begin{equation*}
        \begin{aligned}
            d_{X_2}(\Psi((\ba,\epsilon),t),\Psi((\ba',\epsilon),t))&=d_{X_2}(\Psi((\sigma(\ba),-\epsilon),\widetilde{t}),\Psi((\sigma(\ba'),-\epsilon),\widetilde{t}'))\\
            &\leq 2e^{T_K+1}(|x(\sigma(\ba))-x(\sigma(\ba'))|+|z(\sigma(\ba))-z(\sigma(\ba'))|)+|\widetilde{t}-\widetilde{t}'|\\
            &\leq (2e^{T_K+1}C_1\theta^{-2}+C_1)\theta^{2n}.
        \end{aligned}
    \end{equation*}
    This completes the proof.
\end{proof}

\section{Periodic measures and their arithmetic corrections}\label{S:periodic-arithmetic}

In this section, we establish the periodic approximation of invariant measures for the geodesic flow in $X_2=\SL_2(\ZZ)\backslash \SL_2(\RR)$, and formulate a general proposition that allows us to add certain arithmetic corrections. These arguments are based on the suspension $\Psi_+:\Sigma\times\RR\to X_2$ defined in \eqref{E:+susmap} and the matrix discriminant introduced in Section \ref{SSS:discriminant}.

\subsection{Periodic approximation without arithmetic information} 
We first introduce some notations. Denote by $\mathscr M(X_2)^F$ the set of $F$-invariant probability measures on $X_2$. For a compact subset $K\subset X_2$, put
$
 \sM(K)^F=\{\mu\in\sM(X_2)^F:\supp{\mu}\subset K\}.
$

For an integer $M\geq 2$, we write $$
\Sigma_M:=\{1,\cdots,M\}^{\ZZ},\quad B_M:=\{
(\ba,t):\ba\in\Sigma_M,\;0\leq t\leq R(\ba)
\},\quad \text{and}\quad K_M:=\Psi_+(B_M).
$$
It is clear that each $K_M$ is a compact subset of $X_2$. We first observe the following fact.

\begin{lemma}\label{L:K-correct}
    Let $K$ be a compact subset of $X_2$. Then $\sM(K)^F\subseteq\sM(K_{\exp(T_K)})^F$, where $$T_K=\max\left\{0,\;\log\left(2\max_{q\in K}\het(q)\right)\right\}\quad \text{ is defined in \eqref{def-tk}}.$$
\end{lemma}
\begin{proof}
    Let $\mu\in \sM(K)^F$, that is, $\mu\in\sM(X_2)^F$ and $\supp{\mu}\subseteq K$. It follows that any $q\in \supp(\mu)$ has its full $F$-orbit contained in $K$. We first claim that: $$
    q\in X_2^0 \qquad \text{for }X_2^0\subseteq X_2\text{ defined in \eqref{E:X2-0}}.
    $$ In fact, write $q=\Gamma g$ and first suppose that $\ZZ^2g$ contains a nonzero vector in $\RR e_1$. Then $\lambda_1(\ZZ^2g\cdot a_t)\to 0$ as $t\to -\infty$. By Mahler's compactness criterion, for $t\leq 0$ with $|t|$ large enough, $\Gamma g\cdot a_t\notin K$. This is a contradiction. Similarly, we see that $\ZZ^2g$ cannot contain a nonzero vector in $\RR e_2$. This proves that $q\in X_2^0$.

    Using the suspension \eqref{E:susmap}, we write $q=\Gamma g_{\epsilon}(x(\ba),z(\ba))\cdot a_t$ for some $\ba\in\Sigma$, $\epsilon\in\{\pm 1\}$, and $0\leq t<\tau(x(\ba))$. For $j\in\ZZ$, consider the midpoint of the $j$-th return segment:
\[
q_j:=\Psi\left(
(\sigma^j(\ba),(-1)^j\epsilon),\;\tau(x(\sigma^j(\ba)))/2
\right)\in K.
\]
By Lemma \ref{L:compactdist}, we have $\tau(x(\sigma^j(\ba)))/2\leq T_K$, which implies that $a_{j+1}\leq e^{T_K}$. This proves $\ba\in \Sigma_{\exp(T_K)}$. When $\epsilon=+1$, we see that $q=\Psi_+(\ba,t)\in K_{\exp(T_K)}$. When $\epsilon=-1$, we see that $$
q=\Psi((\ba,-1),t)=\Psi_+(\sigma^{-1}(\ba),\tau(x(\sigma^{-1}(\ba)))+t)\in K_{\exp(T_K)}.
$$
This proves that $\supp{\mu}\subseteq K_{\exp(T_K)}.$
\end{proof}

We remark that the compact set $K$ is not necessarily contained in $K_{\exp(T_K)}$. So we define \begin{equation}\label{E:large-cpt}
    K':=K\cup K_{\exp(T_K)}.
\end{equation}

Now we study the density of all $F$-periodic measures (without arithmetic information) in $\sM(X_2)^F$. 

\begin{lemma}\label{L:no-arith}
    For each $M\geq 2$, the periodic measures $\mu_W$ for
nonempty words $W$ of even length over the finite alphabet $\{1,\ldots,M\}$ are dense in
$\sM(K_M)^F$. Moreover,
\[
 \overline{\bigcup_{M\ge2}\sM(K_M)^F}=\sM(X_2)^F.
\]
The closure here is taken with respect to the weak-${}^*$ topology.
\end{lemma}
\begin{proof}
The suspension $\Psi_+:\Sigma\times \RR\to X_2$ induces, for each $M\geq 2$, a homeomorphism $$
S_M:=B_M/\sim\to K_M,\qquad [(\ba,t)]\mapsto \Psi_+(\ba,t),
$$
where the identification is given by $(\ba,R(\ba))\sim (\sigma^2(\ba),0)$. This homeomorphism conjugates the suspension flow with the diagonal flow, and hence pushes forward any measure on $S_M$ invariant under the suspension flow to a measure on $K_M$ invariant under the diagonal flow. In view of \cite[p.~90]{PP90}, any measure on $S_M$ invariant under the suspension flow is of the form $$
 \widetilde{\nu}=\frac{1}{\int_{\Sigma_M}R\mathrm{d}\nu}\cdot (\nu\times \mathrm{d}t)|_{B_M},
$$
where $\nu$ is a $\sigma^2$-invariant measure on $\Sigma_M$. By \cite[Theorem~3.3]{Par61}, the $\sigma^2$-periodic measures are dense in the set of all $\sigma^2$-invariant measure on $\Sigma_M$. So we conclude that the $F$-periodic measures are dense in the set of all $F$-invariant measure on $K_M$.

Moreover, in view of the above homeomorphism, any $F$-periodic orbit in $K_M$ is of the form $\cO_W$ for some nonempty word $W$ of even length over the finite alphabet $\{1,\cdots,M\}$. This proves the density of all such periodic measures $\mu_W$ in $\sM(K_M)^F$.

Finally, we prove that the union of all $\sM(K_M)^F$ is dense in $\sM(X_2)^F$. Let $\mu\in\sM(X_2)^F$. We first observe that: \begin{equation*}
    \mu(X_2^0)=1\qquad \text{for }X_2^0\subseteq X_2\text{ defined in \eqref{E:X2-0}}.
\end{equation*}
In fact, write $E_1=\{\Gamma g\in X_2:\ZZ^2g\cap\RR e_1\neq \{0\}\}$. Then for any $\Gamma g\in E_1$, $\lambda_1(\ZZ^2g\cdot a_t)\to 0$ as $t\to -\infty$. Let $K\subseteq X_2$ be a compact subset. By Mahler's compactness criterion, for $t\leq 0$ with $|t|$ large enough, $\Gamma g\cdot a_t\notin K$. It follows from the $F$-invariance of $\mu$ and the dominated convergence theorem that $$
\mu(E_1\cap K)=\int_{E_1}1_K(\Gamma g)\mathrm{d}\mu(\Gamma g)=\int_{E_1}1_K(\Gamma g\cdot a_t)\mathrm{d}\mu(\Gamma g)\to 0\quad \text{as }t\to-\infty.
$$
By taking a compact exhaustion of $X_2$, we see that $\mu(E_1)=0$. Similarly, we have $\mu(E_2)=0$ where $E_2=\{\Gamma g\in X_2:\ZZ^2g\cap\RR e_2\neq \{0\}\}$. This proves $\mu(X_2^0)=1$ since $X_2^0=X_2\setminus (E_1\cup E_2)$.

Now consider the restriction of the suspension \eqref{E:+susmap}: $$
\Psi_+\colon\{(\ba,t):\ba\in\Sigma,\;0\leq t<R(\ba)\}\to X_2^0.
$$
%given by
%\begin{equation*}
%\Psi_+(\ba,t)=\begin{cases}
%\Psi((\ba,+),t), &0\le t<\tau(x(\mathbf a)),\\
%\Psi\bigl((\sigma(\ba),-),\,t-\tau(x(\mathbf a))\bigr), &\tau(x(\mathbf a))\le t<R(\mathbf a).
%\end{cases}
%\end{equation*}
It is a Borel isomorphism which conjugates the suspension flow with the diagonal flow. It follows that the $F$-invariant measure $\mu$ restricted on $X_2^0$ can be written as $$
\mu=\frac{1}{\int_{\Sigma}R\mathrm{d}\nu}\cdot \int_{\Sigma}\int_0^{R(\ba)}\delta_{\Psi_+(\ba,t)}\mathrm{d}t\,\mathrm{d}\nu(\ba),
$$
where $\nu$ is a $\sigma^2$-invariant measure on $\Sigma$. By approximating $\nu$ using $\sigma^2$-invariant measures supported on $\Sigma_M\,(M\geq 2)$, we conclude that $\mu$ can be approximated by $F$-invariant measures supported on $K_M\,(M\geq 2)$.
\end{proof}

\subsection{The arithmetic correction of periodic measures} To add some desired arithmetic information on $F$-periodic measures, we carefully modify the tail of each nonempty word with even length to adjust its associated matrix discriminant.

\begin{lemma}\label{L:arith-correct}
    Let $\mathcal P$ be a finite set of primes and $q\notin\mathcal P$ an odd prime. Let $\varepsilon_p\in\{\pm1\}$ be given signs for $p\in\mathcal P$ with $p\ge5$. Then for any nonempty word $W$ of even length, there exists a word $V$ of even length $\leq 2(|G_m|-1)$ where $m=q^2\cdot\prod_{p\in\cP}p$, over the alphabet $\{1,2\}$, such that
\begin{equation}\label{E:force-ramification}
 p\nmid\Delta(WV)\ (p\in\mathcal P),\quad
 v_q(\Delta(WV))=1,\quad
 \text{and}\quad  \chi_p(\Delta(WV))=\varepsilon_p\quad(p\in\mathcal P,\ p\ge5).
\end{equation}
If $2,3\in\mathcal P$, one may also require $d(WV)\equiv5\pmod{24}$.
\end{lemma}
\begin{proof}
We divide the proof into two steps. \medskip

\textsl{Step 1: The congruence correction.} Fix $m\geq 2$ and put $G_m=\SL_2(\ZZ/m\ZZ)$. We show that: for any $g\in G_m$ and any nonempty word $W$ of even length, there exists a word $V$ of even length over the alphabet $\{1,2\}$ such that \begin{equation}\label{E:congru-correct}
|V|\leq 2(|G_m|-1),\quad\text{and}\quad C(W)C(V)\equiv g \pmod m.
\end{equation}
In fact, it is straightforward to see that $G_m$ is generated by the two elements $$
\begin{pmatrix}
    1 &1\\ 0 &1
\end{pmatrix}
=(C(2)C(1))\cdot (C(1)^2)^{-1}\quad \text{and}\quad 
\begin{pmatrix}
    1 &0\\ 1 &1
\end{pmatrix}
=(C(1)^2)^{-1}\cdot (C(1)C(2)).
$$
In particular, the given element $C(W)^{-1}g\in G_m$ can be written as $C(V)$ for some word $V$ of even length over the alphabet $\{1,2\}$. By choosing a shortest path from $1_{G_m}$ to $C(V)$ in $G_m$, we may assume that $|V|\leq 2(|G_m|-1)$. This completes the proof of Step 1.

\medskip 
\textsl{Step 2: The arithmetic correction.} For $p\geq 5$, we choose $z_p\in \FF_p^{\times}\setminus\{1\}$ with $\chi_p(z_p)=\epsilon_p$, and set $$
t_p:=2\frac{1+z_p}{1-z_p}.
$$
It follows that $\chi_p(t_p^2-4)=\epsilon_p$. Take $m=q^2\cdot\prod_{p\in\cP}p$. By Chinese remainder theorem, we may choose $t\mod m$ satisfying $$
\begin{aligned}
t&\equiv t_p\pmod p &&\text{for }p\in\mathcal P,\ p\ge5,\\
t&\equiv0\pmod3 &&\text{if }3\in\mathcal P,\\
t&\equiv1\pmod2 &&\text{if }2\in\mathcal P,\\
t&\equiv2+q\pmod{q^2}.
\end{aligned}
$$
For $g=\begin{pmatrix}
    t &-1\\ 1 &0
\end{pmatrix}\in G_m$, by Step 1, we have chosen a word $V$ of even length over the alphabet $\{1,2\}$ such that \eqref{E:congru-correct} holds. Now we verify that the word $V$ satisfies all desired properties. Note that $$
\Delta(WV)=(\tr C(WV))^2-4\equiv (\tr{g})^2-4=t^2-4 \pmod{m}.
$$

First, for $p\in\cP$ with $p\geq 5$, we have $
\Delta(WV)\equiv t_p^2-4 \pmod{p}$, where $\chi_p(t_p^2-4)=\epsilon_p\neq 0$. If $3\in\cP$, then $\Delta(WV)\equiv 2\pmod{3}$. If $2\in\cP$, then $\Delta(WV)\equiv 1\pmod{2}$. So we conclude that $p\nmid \Delta(WV)$. Second, for the odd prime $q\notin\cP$, we have $\Delta(WV)\equiv 4q\pmod{q^2}$, which implies that $v_q(\Delta(WV))=1$.

Finally, if $2,3\in \cP$, then $\tr C(WV)\equiv t\pmod{m}$ is odd and divisible by $3$. It follows that $$
\Delta(WV)=(\tr C(WV))^2-4\equiv 5\pmod{8},\quad\text{and}\quad \Delta(WV)\equiv 2\pmod{3},
$$
Then the square-free part also satisfies $d(WV)\equiv 5\pmod{8}$ and $d(WV)\equiv 2\pmod{3}$. So Chinese remainder theorem gives $d(WV)\equiv 5\pmod{24}$.
\end{proof}

Now we prove the closing lemma with given modifications on periodic orbits. For a bounded Lipschitz function $f$ on $X_2$, define
\begin{equation*}
 \|f\|_{\mathrm{BL}}:=\|f\|_\infty+\operatorname{Lip}(f),\quad\text{where }\|f\|_{\infty}=\sup_{x\in X_2}|f(x)|,\quad \Lip(f)=\sup_{x\neq x'\in X_2}\frac{|f(x)-f(x')|}{d_{X_2}(x,x')}.
\end{equation*}
For two measures $\mu,\mu'$ on $X_2$, define $$
d_{\BL}(\mu,\mu'):=\sup_{\|f\|_{\mathrm{BL}}\le1}|\mu(f)-\mu'(f)|.
$$
The metric $d_{\BL}$ induces weak convergence on $X_2$ by
\cite[Theorem~11.3.3]{Dudley02}.

\begin{proposition}\label{P:arithmetic-closing}
Let $M\ge2$. There exists some constant $C_M>0$ with the following property. %Let $m\ge2$ and $g\in G_m=\SL_2(\ZZ/m\ZZ)$. 
For any nonempty word $W$ of even length over the alphabet $\{1,\ldots,M\}$, any $k\ge1$, and any word $V$ of even length over the alphabet $\{1,2\}$,
%such that
%\begin{equation}\label{E:closing-convergence}
 %|V_n|\le 2(|G_m|-1),\qquad C(W)^nC(V_n)\equiv g\pmod m.
%\end{equation}
%Moreover, there exists some constant $C_M>0$, depending only on $M$, such that
\begin{equation}\label{E:closing-rate}
 d_{\BL}(\mu_{W^kV},\mu_{W})\le\frac{C_M(1+|V|)}{k|W|}.
\end{equation}
\end{proposition}
\begin{proof}
%For $g\in G_m$, $n\geq 1$, and the word $W^n$, let $V_n$ be the word given by Step 1 in the proof of Lemma \ref{L:arith-correct}. We verify that the new word $W^nV_n$ satisfies the estimate \eqref{E:closing-rate}.
Let $n_0>0$ be the absolute constant given by Proposition \ref{P:cpt-shadow}. We first show that: there exists some constant $C_M^0>0$ such that for $\ba,\ba'\in\Sigma_M$ whose digit numbers agree on $[-h,h]$ where $h\geq n_0+1$, and any bounded Lipschitz function $f$ on $X_2$, one has \begin{equation}\label{E:Jf-close}
    |J_f(\ba)-J_f(\ba')|\leq C_M^0\|f\|_{\BL}\theta^{2h}.
\end{equation}
In fact, for $i=0,1$, write $\tau_i=\min\{\tau(x(\sigma^i(\ba))),\;\tau(x(\sigma^i(\ba')))\}$. Clearly, $|\tau_i|\leq 2\log(M+1)$. By Proposition \ref{P:cpt-shadow}, for any $0\leq t<\tau_i$, we have $$
d_{X_2}(\Psi((\sigma^i(\ba),(-1)^i),t)),\;\Psi((\sigma^i(\ba'),(-1)^i),t))\leq C_{K_M}\theta^{2(h-i)}.
$$
It follows that \begin{equation*}
        \left|\int_0^{\tau_i}f(\Psi((\sigma^i(\ba),(-1)^i),t))\mathrm{d}t-\int_0^{\tau_i}f(\Psi((\sigma^i(\ba'),(-1)^i),t))\mathrm{d}t\right|\leq \tau_i\cdot \Lip(f)C_{K_M}\theta^{2(h-i)}.
\end{equation*}
By Corollary \ref{C:samepartial}, the difference of integrals on other parts is at most $$
|\tau(x(\sigma^i(\ba)))-\tau(x(\sigma^i(\ba')))|\cdot \|f\|_{\infty}\leq C_1\|f\|_{\infty}\theta^{2(h-i)}.
$$
Therefore, we conclude that $$
|J_f(\ba)-J_f(\ba')|\leq \max\{2\log(M+1)C_{K_M},\,C_1\}(1+\theta^{-2})\cdot \|f\|_{\BL}\theta^{2h}.
$$

Next, we show that: there exists some constant $C_M^1>0$ such that for any $k\geq 1$ and any bounded Lipschitz function $f$ on $X_2$, one has \begin{equation}\label{E:If-close}
    |I_f(W^kV)-k\cdot I_f(W)|\leq (C_M^1+2\log{3}\cdot |V|)\|f\|_{\BL}.
\end{equation}
In fact, write $2\ell=k|W|$ and $2s=|V_n|$. Then we have $$
I_f(W^kV)-k\cdot I_f(W)=\sum_{j=0}^{\ell-1}\left(J_f(\sigma^{2j}(\overline{W^kV}))-J_f(\sigma^{2j}(\overline{W}))\right)+\sum_{j=\ell}^{\ell+s-1}J_f(\sigma^{2j}(\overline{W^kV})).
$$
For $0\leq j\leq \ell-1$, note that the digit numbers of the words $\overline{W^kV}$ and $\overline{W}$ agree on $[-h_j,h_j]$, where $h_j=\min\{2j-1,\,2\ell-2j\}$. When $h_j<n_0+1$, one has $$
\left|J_f(\sigma^{2j}(\overline{W^kV}))-J_f(\sigma^{2j}(\overline{W}))\right|\leq 2\cdot 4\log(M+1)\cdot \|f\|_{\infty}.
$$
When $h_j\geq n_0+1$, by \eqref{E:Jf-close}, one has $$
\left|J_f(\sigma^{2j}(\overline{W^kV}))-J_f(\sigma^{2j}(\overline{W}))\right|\leq C_M^0\|f\|_{\BL}\theta^{2h_j}.
$$
So we have \begin{equation*}
    \begin{aligned}
        \sum_{j=0}^{\ell-1}\left|J_f(\sigma^{2j}(\overline{W^kV}))-J_f(\sigma^{2j}(\overline{W}))\right|
        &\leq 8(n_0+1)\log(M+1)\|f\|_{\infty}
        +\sum_{\substack{0\leq j\leq \ell-1\\h_j\geq n_0+1}}C_M^0\|f\|_{\BL}\theta^{2h_j}\\
        &\leq \left(8(n_0+1)\log(M+1)+2\frac{\theta^{2(n_0+1)}}{1-\theta^2}C_M^0\right)\|f\|_{\BL}.
    \end{aligned}
\end{equation*}
On the other hand, since $V$ has alphabet $\{1,2\}$, we see that for any $0\leq j\leq s-1$, $R(\sigma^{2j}(\overline{VW^k}))\leq 4\log{3}$. Then we have $$
\left|\sum_{j=\ell}^{\ell+s-1}J_f(\sigma^{2j}(\overline{W^kV}))\right|=\left|\sum_{j=0}^{s-1}J_f(\sigma^{2j}(\overline{VW^k}))\right|\leq 4\log{3}\cdot s\|f\|_{\infty}.
$$
Putting these two inequalities together gives \eqref{E:If-close}.

Finally, for any bounded Lipschitz function $f$ on $X_2$ with $\|f\|_{\BL}\leq 1$, one has \begin{equation*}
    \begin{aligned}
        |\mu_{W^kV}(f)-\mu_W(f)|
        &=\frac{1}{T(W^kV)}\left|
I_f(W^kV)-k\cdot I_f(W)-(T(W^kV)-k\cdot T(W))\cdot \mu_W(f)
\right|\\
&\leq \frac{2}{|W^kV|\cdot \log{2}}(C^1_M+2\log{3}\cdot |V|).
    \end{aligned}
\end{equation*}
This is as desired.
\end{proof}

\section{Density of periodic measures with arithmetic information}\label{S:arithmetic-density}

In this section, we prove a stronger result concerning the density of periodic measures with arithmetic information, from which Theorem \ref{T:main1} immediately follows. 

\begin{theorem}\label{T:main-arithmetic}
Let $A\subset\mathcal S_f$ be a finite set, let $h_1,h_2:\mathbb N\to\mathbb N$ satisfy $h_1(n),h_2(n)\to +\infty$ as $n\to\infty$, and fix signs $\varepsilon_p\in\{\pm1\}$ for all primes $p\ge5$. Then there exist nonempty words $U_n\,(n\geq 1)$ of even lengths with the following properties. Write
\[
 \Delta_n=\Delta(U_n),\qquad d_n=d(U_n).
\]
We have
\begin{enumerate}
\item $B=\{d_n:n\geq 1\}$ is a strictly increasing sequence in $\cS_f\setminus A$, and
\[
 \lim_{T\to +\infty}\frac{|B\cap[0,T]|}{h_1(T)}=0.
\]
\item for all $n\ne n'$, $\gcd(\Delta_n,\Delta_{n'})=1$, and every prime divisor $p$ of $\Delta_n$ satisfies $p>h_2(n)$. 
\item for all $n\geq 1$, $d_n\equiv5\pmod{24}$; for all prime $p\ge5$ and all sufficiently large $n$, $\chi_p(d_n)=\epsilon_p$.
Thus $2$ and $3$ are inert in all $\mathbb Q(\sqrt{d_n})$, while any prime $p\ge5$ is split or inert in $\mathbb Q(\sqrt{d_n})$ for all sufficiently large $n$, according to $\varepsilon_p=1$ or $-1$.
\item for any $M\ge2$ and $N\ge1$,
\[
 \sM(K_M)^F\subset
 \overline{\{\mu_{U_n}:n\ge N,\;\supp\mu_{U_n}\subset K_M\}}.
\]
%In particular, each tail of $(\rho_n)$ is dense in $\mathscr M(X_2)$.
\item for any compact $K\subset X_2$ and any $\mu\in\mathscr M(K)^F$, there is a subsequence $\mu_{U_{n_j}}\to\mu$ with
$\supp\mu_{U_{n_j}}\subset K'$ for all $j\geq 1$, where $K'$ is defined in \eqref{E:large-cpt}.
\end{enumerate}
\end{theorem}

\begin{proof}
Without loss of generality, we assume that $h_1,h_2$ are both increasing functions. Set $A_0=\max(A\cup\{1\})$. For $n\geq 1$, choose $R_n\in\NN$ with $h_1(R_n)\ge n^2$. We enumerate the set $$
\{(M,W):M\geq 2,\;W\text{ is a nonempty word of even length over the alphabet }\{1,\ldots,M\}\}
$$
as a sequence $\{(M_n,W_n):n\geq 1\}$, in which every pair appears infinitely often.

We proceed recursively. Put $d_0=1$. For $n\geq 1$, suppose that $U_1,\ldots,U_{n-1}$ are chosen. Let
\[
 \mathcal P_n=\{2,3\}\cup\left\{p\text{ prime}:p\le h_2(n)\right\}
 \cup\bigcup_{j\leq n-1}\{p\text{ prime}:p\mid\Delta_j\}.
\]
We also choose a prime $q_n\notin\mathcal P_n$ with $q_n>\max\{A_0,R_n,d_{n-1}\}$. Let $m_n=q_n^2\cdot\prod_{p\in\cP_n}p$ and let $C_{M_n}>0$ be given by Proposition \ref{P:arithmetic-closing}. Choose $k_n\geq 1$ sufficiently large so that \[
 \frac{C_{M_n}(1+2(|G_{m_n}|-1))}{k_n|W_n|}<\frac{1}{n}.
\]
Then by Lemma \ref{L:arith-correct}, there exists a word $V_n$ of even length $\leq 2(|G_{m_n}|-1)$ over the alphabet $\{1,2\}$, such that for $U_n:=W_n^{k_n}V_n$, $$
p\nmid\Delta_n\ (p\in\mathcal P_n),\quad
 v_{q_n}(\Delta_n)=1,\quad
 \chi_p(\Delta_n)=\varepsilon_p\quad(p\in\mathcal P_n,\ p\ge5),\quad 
 \text{and}\quad d_n\equiv 5\pmod{24}.
$$
It is clear that $\mu_{U_n}\in\sM(K_{M_n})^F$ and $q_n\mid d_n$. By Proposition \ref{P:arithmetic-closing}, we have
\[
 d_{\BL}(\mu_{U_n},\mu_{W_n})\leq \frac{C_{M_n}(1+|V_n|)}{k_n|W_n|}<\frac1n.
\]
These prove (2) and (3).  

Since $d_n\geq q_n>\max\{A_0,R_n,d_{n-1}\}$, we see that $B=\{d_n:n\geq 1\}$ is a strictly increasing sequence in $\cS_f\setminus A$. Moreover, for $n_T=|B\cap [0,T]|$, we have $T\geq d_{n_T}>R_{n_T}$, and hence $$
\frac{|B\cap[0,T]|}{h_1(T)}\leq \frac{n_T}{h_1(R_{n_T})}\leq \frac{1}{n_T}\to 0\quad\text{as } T\to +\infty.
$$
These prove (1). 

Fix $M\geq 2$ and $N\geq 1$. For each nonempty word $W$ of even length over the alphabet $\{1,\ldots,M\}$, there are arbitrarily large $n\geq N$ with $(M_n,W_n)=(M,W)$, for which $\mu_{U_n}$ are supported in $K_{M_n}=K_M$ and converge to $\mu_{W_n}=\mu_W$. This proves the density of $$
\{\mu_{U_n}:n\ge N,\;\supp\mu_{U_n}\subset K_M\}
$$
in the following set consisting of all periodic measures supported in $K_M$:
$$\{\mu_W:W\text{ is a nonempty word of even length over the alphabet }\{1,\cdots,M\}\}.$$
By Lemma \ref{L:no-arith}, the latter set is dense in $\sM(K_M)^F$. This proves (4).

For any compact set $K$, we see from Lemma \ref{L:K-correct} that $\sM(K)^F\subseteq \sM(K_{\exp(T_K)})^F$. Since $K'=K\cup K_{\exp(T_K)}$, the conclusion (5) follows from (4).   
\end{proof}

\begin{remark}\label{R:strong}
    Theorem \ref{T:main-arithmetic} strengthens Theorem \ref{T:main1} in the following two aspects: \begin{enumerate}
        \item[(i)] We may require that the matrix discriminants $\Delta_n$ and their square-free parts $d_n$ satisfy (2) and (3), with a prescribed function $h_2:\NN\to\NN$ and associated signs $\epsilon_p\in\{\pm 1\}$ for all primes $p\geq 5$;
        \item[(ii)] We find a uniform sequence of invariant measures $\{\mu_{U_n}:n\geq 1\}$, such that for any compact subset $K\subseteq X_2$, any invariant measure supported in $K$ can be approximated by a subsequence $\{\mu_{U_{n_j}}:j\geq 1\}$ whose supports are contained in a larger compact set (only determined by $K$).
    \end{enumerate}
\end{remark}

\section{Further discussion}\label{S:discuss}

In this section, we propose higher-rank analogues of Theorem \ref{T:main1} and Problem \ref{P:fixed-field-density}, which provided the initial impetus for this paper. Let $n\ge3$, and set
\[
 G=\SL_n(\RR),\qquad \Gamma=\SL_n(\ZZ),\qquad
 X_n=\Gamma\backslash G.
\]
Here $G$ acts on $X_n$ by right multiplication. Let $A\leq G$ denote the connected positive diagonal subgroup:
\[
 A:=\left\{\diag(a_1,\ldots,a_n):a_i>0,\ \prod_{i=1}^n a_i=1
 \right\}.
\]
Let $F=\{g_t:t\in\RR\}$ be a one-parameter subgroup in $A$, and write $g_t=\diag(e^{t\lambda_1},\ldots,e^{t\lambda_n})$ where $\sum_{i=1}^n\lambda_i=0.$ We assume that $F$ is \textsl{regular}, meaning that $\lambda_i\ne\lambda_j$ for every $i\ne j$ and that $\lambda_i\ne0$ for every $i$. Write \begin{subequations}
    \begin{align*}
    \mathscr M(X_n)_c^F&:=\{F\text{-periodic measures on $X_n$}\};\\
        \mathscr M(X_n)^F&:=\{F\text{-invariant measures on $X_n$}\}.
    \end{align*}
\end{subequations}
As in Section \ref{S:intro}, we raise the following problem: 

\begin{problem}\label{P:higher-rank-density}
Assume that $\mathscr M(X_n)_c^F\neq \varnothing$. Is $\mathscr M(X_n)_c^F$ dense in $\mathscr M(X_n)^F$?
\end{problem}

We remark that one of the main motivations for this problem comes from the work of Marcus \cite{Mar}, where the density of periodic measures in all invariant measures was proved for ergodic (not necessarily hyperbolic) toral automorphisms.

\subsection{Periodic points and compact $A$-orbits}

We present several important consequences of the existence of one $F$-periodic point. The first one is a topological result based on the density of commensurator.

\begin{lemma}\label{L:commensurator-periodic}
Let $G$ be a Lie group and $\Gamma<G$ a discrete subgroup. Suppose that the commensurator
\[
 \Comm_G(\Gamma):=\left\{h\in G:\Gamma\cap h^{-1}\Gamma h\text{ has finite index in both }\Gamma\text{ and }h^{-1}\Gamma h\right\}
\]
is dense in $G$. Let $F=\{g_t:t\in\RR\}\leq G$ be a one-parameter subgroup acting on $\Gamma\backslash G$ by right multiplication. If $F$ has a periodic point, then its periodic points are dense in $\Gamma\backslash G$.
\end{lemma}

\begin{proof}
Let $x=\Gamma g$ be an $F$-periodic point. Choose $T>0$ so that $xg_T=x$, and set
\[ \gamma:=g g_T g^{-1}\in\Gamma.\]
For any $h\in\Comm_G(\Gamma)$, the subgroup $\Gamma\cap h^{-1}\Gamma h$ has finite index in $\Gamma$. Then there exists $k\geq1$ such that $\gamma^k\in\Gamma\cap h^{-1}\Gamma h.$ Thus we have \[
 (hg)g_{kT}(hg)^{-1}=h\gamma^kh^{-1}\in\Gamma.
\]
It follows that $\Gamma hg$ is an $F$-periodic point. Since $\Comm_G(\Gamma)$ is dense in $G$, the points $\Gamma hg$ obtained in this way are dense in $\Gamma\backslash G$. The proof is complete.
\end{proof}

\begin{remark}
In our setting $G=\SL_n(\RR),\Gamma=\SL_n(\ZZ),X_n=\Gamma\backslash G,$
the commensurator $\Comm_G(\Gamma)$ contains a dense subgroup $\SL_n(\QQ)$ of $G$. Hence, if a one-parameter subgroup $F\leq G$ has a periodic point in $X_n$, then its periodic points are dense in $X_n$.
\end{remark}

The second one is an arithmetic consequence based on the classification of compact $A$-orbits. Here we introduce a concept stronger than regularity.

\begin{definition}\label{D:good-regular}
A regular one-parameter subgroup $F=\{g_t:t\in\RR\}$ in $A$ is good regular, if for $g_t=\diag(e^{t\lambda_1},\cdots,e^{t\lambda_n})$, one has
\[
 \sum_{i\in J}\lambda_i\ne0\qquad
 \text{ for any }\varnothing\ne J\subsetneqq\{1,\ldots,n\}.
\]
Equivalently, for any $t\ne0$, no nonempty proper subset of the eigenvalues of $g_t$ has product $1$.
\end{definition}

\begin{lemma}\label{L:good-regular-compact-orbit}
Let $F\leq A$ be a regular one-parameter subgroup. Then 
\begin{enumerate}
 \item For every $L>0$, there are only finitely many compact $A$-orbits in $X_n$ that contain an $F$-periodic point of primitive period at most $L$.
\item If $F$ is good regular, then every $F$-periodic point in $X_n$ belongs to a compact $A$-orbit. Consequently, for every $L>0$, the $F$-periodic points of primitive period at most $L$ are contained in a finite union of compact $A$-orbits.
\end{enumerate}
\end{lemma}
\begin{proof}
(1) Let $xA$ be a compact $A$-orbit containing an $F$-periodic point $x$ with return time $0<T\leq L$. By the arithmetic correspondence for compact $A$-orbits (see \cite[Theorem~3.1]{Sha11}), $xA$ has an associated totally real field $K$ of degree $n$ and an associated order $\cO\subseteq K$. The return element $g_T$ corresponds to a unit $\alpha\in\cO$ whose real conjugates are $e^{T\lambda_1},\ldots,e^{T\lambda_n}.$ Since $F$ is regular, these conjugates are distinct. Hence $K=\QQ(\alpha)$ and $\ZZ[\alpha]\subseteq\cO$ is an order in $K$. It follows that
\[
 \disc(\cO)\leq \disc(\ZZ[\alpha])=\prod_{i<j}\bigl(e^{T\lambda_i}-e^{T\lambda_j}\bigr)^2\ll_{F,L}1.
\]
By \cite[Proposition~2.9 and Corollary~4.4]{ELMV09}, there are only finitely many compact $A$-orbits in $X_n$ whose associated orders have bounded discriminant. This proves (1).

\smallskip
\noindent
(2) Suppose now that $F$ is good regular, and let $x=\Gamma g$ be an $F$-periodic point. Choose $T>0$ such that $xg_T=x$, and put\[
 \gamma:=gg_Tg^{-1}\in\Gamma.
\]
We claim that: its characteristic polynomial $P_\gamma(z):=\det(zI-\gamma)\in\ZZ[z]$ is irreducible over $\QQ$.

Suppose to the contrary that $P_\gamma(z)$ is not irreducible over $\QQ$. By Gauss' lemma, $P_\gamma$ has a monic factor $P\in\ZZ[z]$ with $0<\deg{P}<\deg{P_{\gamma}}$ and constant term $\pm1$. Since the roots of $P_\gamma$ are $e^{T\lambda_1},\cdots,e^{T\lambda_n}$, we see that the roots of $P$ correspond to a nonempty proper subset $J\subsetneqq\{1,\cdots,n\}$, and their product is
\[
 \prod_{j\in J}e^{T\lambda_j}
 =(-1)^{\deg P}\cdot P(0)=1.
\]
It follows that $\sum_{j\in J}\lambda_j=0$, which is impossible due to the good regularity of $F$. Therefore, $P_\gamma$ is irreducible over $\QQ$.

By the arithmetic correspondence for compact $A$-orbits (see \cite[Theorem~3.1]{Sha11}), the orbit $xA$ is compact. The final assertion now follows from (1).
\end{proof}

For good regular $F$, the above lemma places every $F$-periodic orbit inside a compact $A$-orbit. This naturally leads to a higher-rank fixed-field question for the full $A$-action on $X_n$.

\begin{problem}\label{P:fixed-field-A-density}
Fix a totally real field $K$ of degree $n$. Is the set of normalized invariant measures on compact $A$-orbits associated with $K$ dense in the set of all $A$-invariant measures on $X_n$?
\end{problem}

\begin{remark}\label{R:fixed-field-A-Haar}
The Haar probability measure on $X_n$ is known to be approximable by individual compact $A$-orbit measures associated with any fixed totally real field $K$ of degree $n$. This follows from Shapira and Zheng's equidistribution theorem \cite[Theorem~1.2]{SZ24}. %Thus the fixed-field problem asks whether such approximation is possible for every $A$-invariant measure, rather than only for $m_n$.

The same conclusion can also be obtained by a Hecke argument: the Hecke transforms of a compact $A$-orbit associated with $K$ decompose into finitely many compact $A$-orbits associated with the same field. By Hecke equidistribution \cite[Theorem~1.2]{EO06} and the ergodicity of the Haar measure, one can select a sequence of these individual orbit measures converging to the Haar measure on $X_n$.
\end{remark}

For related work on limits of measures supported on compact
$A$-orbits, including non-ergodic limits and escape of mass in the vague topology, see Solan and Yifrach \cite{SY23}. In particular, \cite[Open Question~1.12]{SY23} asks whether a compact $A$-orbit measure can be approximated by measures supported on other compact $A$-orbits.

\subsection{Return times for a good regular diagonal flow}
In this subsection, we prove a finer structural result on the set of return times for a good regular diagonal flow. They are closely related to the famous six exponentials theorem \cite[Theorem 10]{Wald08} and the four exponentials conjecture \cite[Conjecture 12]{Wald08}.

Let $F=\{g_t:t\in\RR\}$ be a one-parameter subgroup in $A$, and write $$
\Per(F):=\{t\in\RR: xg_t=x\text{ for some }x\in X_n\}.
$$
This topological information can also be recorded in the following number-theoretic way. For $t\in\RR$, write $P_t(z):=\det(zI_n-g_t)=\prod_{j=1}^n(z-e^{t\lambda_j})$. Set 
 $$\cT(F):=\{t\in\RR:P_t(z)\in\ZZ[z]\text{ and }P_t(z)\text{ is irreducible}\}.
$$
It is clear that for $t\in\cT(F)$, $K_t:=\QQ(e^{t\lambda_1})$ is totally real field of degree $n$. Set $$\cK(F):=\{K_t:t\in\cT(F)\}.$$
We start with a number-theoretic description of $\Per(F)$. 

\begin{lemma}\label{L:field-time-unit}
Let $F\leq A$ be a good regular one-parameter subgroup. 
\begin{enumerate}
\item We have $\Per(F)=\{0\}\cup\cT(F)$ is discrete;
 \item For $s,t\in\cT(F)$, if $s/t\in\QQ^{\times}$, then $K_s=K_t$;
 \item For every $t_0\in\cT(F)$, write $$
 \cT_{t_0}(F):=\{0\}\cup\bigl(\cT(F)\cap\QQ t_0\bigr).
 $$
 Then there exists $\tau_{t_0}>0$ such that $\cT_{t_0}(F)=\tau_{t_0}\ZZ.$
\end{enumerate}
\end{lemma}
\begin{proof}
(1) This is the arithmetic correspondence for compact $A$-orbits (see \cite[Theorem~3.1]{Sha11}).
~\\
    (2) For $t\in\cT(F)$, we see that $e^{t\lambda_1}$ is a totally positive algebraic number. It follows that any complex embedding of the field $\QQ(e^{t\lambda_1})$ is determined by any of its positive power. So we have $\QQ(e^{t\lambda_1})=\QQ(e^{qt\lambda_1})$ for any $q\in\ZZ\setminus\{0\}$, and hence $\QQ(e^{t\lambda_1})=\QQ(e^{s\lambda_1})$ for any $s/t\in\QQ^{\times}$.
~\\
    (3) We first verify that $\cT_{t_0}(F)$ is a subgroup of $(\RR,+)$. In fact, let $s,t\in \cT_{t_0}(F)$ and suppose that $s,t,s+t\neq 0$. It follows from $s,t\in\QQ t_0$ that $s+t\in\QQ t_0$ and that $K:=K_s=K_t=K_{s+t}$. Moreover, since $s,t\in\cT(F)$, we see that the element $e^{(s+t)\lambda_1}=e^{s\lambda_1}\cdot e^{t\lambda_1}\in\cO_K^{\times}$ has $n$ distinct conjugates $e^{(s+t)\lambda_j}\,(1\leq j\leq n)$. This proves that $s+t\in \cT(F)$. A similar argument show that if $t\in\cT_{t_0}(F)$, then $-t\in \cT_{t_0}(F)$. This verifies the subgroup axioms.

    Finally, we verify that the set $\{0\}\cup\cT(F)$ is discrete in $\RR$. In fact, since the function $$
    \RR\to\RR_{>0},\qquad t\mapsto \sum_{j=1}^ne^{t\lambda_j}.
    $$
    is strictly convex and takes values in $\ZZ$ when restricted to $\cT(F)$, we see that the set $\cT(F)$ is locally finite. Thus, the set $\{0\}\cup\cT(F)$ is discrete. In particular, its subset $\cT_{t_0}(F)$ is also discrete. We conclude that $\cT_{t_0}(F)=\tau_{t_0}\ZZ$ for some $\tau_{t_0}>0$.
\end{proof}

\begin{proposition}\label{P:periodic-points}
    Let $F\leq A$ be a good regular one-parameter subgroup. Suppose that $F$ has a periodic point. \begin{itemize}
 \item[(1)] If $n\geq4$, then $\Per(F)=\tau\ZZ$ for some $\tau>0$; moreover, $|\cK(F)|=1$;
 \item[(2)] If $n=3$, then either $\Per(F)=\tau\ZZ$ for some $\tau>0$, or $\Per(F)=\tau_1\ZZ\cup\tau_2\ZZ$ for some $\tau_1,\tau_2>0$ with $\tau_1/\tau_2\notin\QQ$; moreover, $|\cK(F)|\leq 2$.
 \end{itemize}   
\end{proposition}

We remark here that the item (1) is a consequence of the six exponentials theorem, while in the item (2), if the four exponentials conjecture holds, then we only have the first possibility.

\begin{proof}[Proof of Proposition \ref{P:periodic-points}] 
We first claim that: for $1\leq i\neq j\leq n$, the real numbers $\lambda_i$ and $\lambda_j$ are linearly independent over $\QQ$. In fact, let $t\in\cT(F)$. Suppose that $q\lambda_i=p\lambda_j$ for some $i\neq j$ and $p,q\in\ZZ\setminus\{0\}$. Then we have $(e^{t\lambda_i})^q=(e^{t\lambda_j})^p$. Taking their positive conjugates gives $$
\{(e^{t\lambda_k})^q:1\leq k\leq n\}=\{(e^{t\lambda_\ell})^p:1\leq \ell\leq n\}.
$$
In particular, by taking logarithms and comparing
the sums of their squares, we obtain $q=\pm p$ and hence $\lambda_i=\pm \lambda_j$, which contradicts the good regularity. This proves our claim.

Write $V=\Span_{\QQ}(\{\lambda_1,\cdots,\lambda_n\})$ and $r=\dim_{\QQ}{V}.$ The above claim gives $r\geq 2$.

\medskip\noindent 
(1) Let $n\geq 4$. We claim that: $r\geq 3$. Suppose that $r=2$. Let $t\in\cT(F)$ and let $S$ be the Galois group of the splitting field of $P_t(z)$ over $\QQ$. Then the group action induces a faithful representation $S\to \GL(tV)$ that transitively permutes $t\lambda_1,\cdots,t\lambda_n$. In particular, we have $|S|\geq n\geq 4$. Since every element of $S$ has determinant $\pm 1$, we see that the kernel $\{\sigma\in S:\det{\sigma}=1\}$ is nontrivial. By Cauchy's theorem, this kernel contains an element $\sigma_0$ of prime order $p$. It follows that $\sigma_0\in\GL(tV)$ has a primitive $p$-th root of unity as an eigenvalue. Since $\dim_{\QQ}{V}=2$, we see that $p\leq 3$.

If $p=2$, then $\id_{tV}+\sigma_0=0$. In particular, we have $t\lambda_i+t\lambda_j=t\lambda_i+\sigma_0(t\lambda_i)=0$ for some $i,j$, which contradicts the good regularity. If $p\geq 3$, then $\id_{tV}+\sigma_0+\sigma_0^2=0$. In particular, we have $t\lambda_i+t\lambda_j+t\lambda_k=t\lambda_i+\sigma_0(t\lambda_i)+\sigma_0^2(\lambda_i)=0$ for some $i,j,k$, which also contradicts the good regularity. This verifies the claim.

Now let $s,t\in\cT(F)$ and take $\lambda_i,\lambda_j,\lambda_k$ to be linearly independent over $\QQ$. Since all these six numbers $e^{s\lambda_i},e^{s\lambda_j},e^{s\lambda_k},e^{t\lambda_i},e^{t\lambda_j},e^{t\lambda_k}$ are algebraic numbers, by the six exponentials theorem \cite[Theorem~10]{Wald08}, we see that $s,t$ must be linearly independent over $\QQ$. This prove that for any $t_0\in\cT(F)$, $\cT_{t_0}(F)=\{0\}\cup\cT(F)$. So the conclusion follows from Lemma \ref{L:field-time-unit}.

\medskip\noindent 
(2) Let $n=3$. Since $r\geq 2$ and $\lambda_1+\lambda_2+\lambda_3=0$, we have $r=2$. Applying the six exponentials theorem to $\lambda_1,\lambda_2$ and any three numbers in $\cT(F)$, we see that $\dim_{\QQ}\Span_{\QQ}\cT(F)\le2.$
Moreover, we claim that: \begin{equation}\label{E:line-union}
\cT(F)\text{ is contained in the union of at most two rational lines.}
\end{equation}

Assume the truth of this claim for one moment. If $\cT(F)\subseteq \QQ t_0$ for some $t_0\in \cT(F)$, then by Lemma \ref{L:field-time-unit} (1) and (3), we have $\Per(F)=\{0\}\cup\cT(F)=\cT_{t_0}(F)=\tau_{t_0}\ZZ$ for some $\tau_{t_0}>0$; moreover, we see that $|\cK(F)|=1$ by Lemma \ref{L:field-time-unit} (2). If $\cT(F)\subseteq \QQ t_0\cup\QQ t_1$ for some $t_0,t_1\in\cT(F)$ with $t_0/t_1\notin\QQ$, then by Lemma \ref{L:field-time-unit} (1) and (3), we have $$
\Per(F)=\{0\}\cup\cT(F)=\cT_{t_0}(F)\cup\cT_{t_1}(F)=\tau_{t_0}\ZZ\cup\tau_{t_1}\ZZ\quad \text{ for some }\tau_{t_0},\tau_{t_1}>0;
$$
moreover, we see that $|\cK(F)|\leq 2$ by Lemma \ref{L:field-time-unit} (2). This will complete the proof.

In the following, we verify the above claim \eqref{E:line-union}. Suppose, to the contrary that, there exist $t_0,t_1,t_2\in\cT(F)$ with $t_1/t_2\notin\QQ$ and $c_1,c_2\in\QQ^{\times}$ such that $t_0=c_1t_1+c_2t_2$. Let $S$ be the Galois group of the splitting field of the three polynomials $P_{t_0},P_{t_1},P_{t_2}$. For $\sigma\in S$ and $i\in\{0,1,2\}$, define
$\pi_i(\sigma)\in S_3$ by
\begin{equation}\label{E:def-pi}
 \sigma(e^{t_i\lambda_j})=e^{t_i\lambda_{\pi_i(\sigma)(j)}}, \qquad\text{for all } 1\leq j\leq 3.
\end{equation}
We show that for any $\sigma\in S$,
\begin{equation}\label{E:cubic-same-permutations}
 \pi_0(\sigma)=\pi_1(\sigma)=\pi_2(\sigma).
\end{equation}

In fact, we see from the identity $t_0=c_1t_1+bc_2t_2$ and \eqref{E:def-pi} that \[
 t_0\lambda_{\pi_0(\sigma)(j)}
 =c_1 t_1\lambda_{\pi_1(\sigma)(j)}+c_2 t_2\lambda_{\pi_2(\sigma)(j)},\qquad\text{for all } 1\leq j\leq 3.
\]
For $i=0,1,2$, write $
 \Lambda_i:=\bigl(\lambda_{\pi_i(\sigma)(1)},
       \lambda_{\pi_i(\sigma)(2)},
       \lambda_{\pi_i(\sigma)(3)}\bigr).$
Then $\|\Lambda_0\|=\|\Lambda_1\|=\|\Lambda_2\|$, and $(c_1t_1+c_2t_2)\Lambda_0=at_1\Lambda_1+bt_2\Lambda_2$. It follows from taking norms that $c_1t_1\cdot c_2t_2\cdot\|\Lambda_1-\Lambda_2\|^2=0$. Thus we conclude that $\Lambda_0=\Lambda_1=\Lambda_2$. This proves \eqref{E:cubic-same-permutations}.

It follows from \eqref{E:cubic-same-permutations} that for any $m_1,m_2\in\ZZ$ and $u:=m_1t_1+m_2t_2$, the Galois conjugates of $$
e^{u\lambda_1}=(e^{t_1\lambda_1})^{m_1}\cdot (e^{t_2\lambda_1})^{m_2}
$$
are $\{e^{u\lambda_1},e^{u\lambda_2},e^{u\lambda_3}\}$. This means that $u\in\{0\}\cup \cT(F)$. So we have $\ZZ t_1+\ZZ t_2\subseteq \{0\}\cup\cT(F)$. However, since $t_1/t_2\notin\QQ$, we may find nonzero elements in $\ZZ t_1+\ZZ t_2$ arbitrarily close to the origin. This contradicts Lemma \ref{L:field-time-unit} (1). Therefore, the claim \eqref{E:line-union} holds.
\end{proof}

%\blue{I think theorem 6.9, 6.10, 6.12 can be greatly simplified, but since you are always much better on editing... I will not touch...}--Yes, I have done so

%\begin{remark}
%We learned the materials concerning four exponential conjecture and six exponential theorem from Antoine Marnat, although the proof here are based on conversations with AI.  (I think it is necessary to acknowledge him, since I really have discussed this problem with him)
%\end{remark}--Yes, I've added this part to acknowledgment

%\enlargethispage{3\baselineskip}

\end{document}